\documentclass[12pt]{amsart}
\usepackage{amsfonts}
\usepackage{amsmath}
\usepackage{amssymb}
\usepackage{graphicx}
\usepackage{hyperref}

\theoremstyle{plain}

\theoremstyle{definition}

\theoremstyle{remark}

\numberwithin{equation}{section}

\renewcommand\epsilon\varepsilon
\renewcommand\phi\varphi

\title[Three-spheres theorem for harmonic functions]{Three-spheres theorem for harmonic functions\\(non-concentric case)}

\author[N.~U.~Arakelian]{Norair U.~Arakelian}
\address{Institute of Mathematics, National Academy of Sciences of Armenia, Yerevan, Armenia}

\author[N.~Matevosyan]{Norayr Matevosyan}
\address{Independent Researcher}
\email{norayr@gmail.com}
\thanks{N.~Matevosyan: ORCID \href{https://orcid.org/0009-0001-2994-9522}{0009-0001-2994-9522}.}

\subjclass[2020]{Primary 31B05; Secondary 31B25, 35B60}
\keywords{Three-spheres theorem, harmonic functions, non-concentric spheres, Hadamard, log-convexity, propagation of smallness}

\date{}

\begin{document}
\begin{abstract}
A direct analog of Hadamard's three-circle theorem is obtained for harmonic functions
(in weighted $L^{2}$-norm) in case of $(n-1)$-dimensional non-concentric spheres in $\mathbb{R}^{n}$.
The result extends the concentric case to correlated non-concentric, non-touching spheres
via an inversion technique. Applications to propagation of smallness and uniqueness
for harmonic functions are given.
\end{abstract}
\maketitle

\renewcommand{\thefootnote}{\fnsymbol{footnote}}
\footnotetext[1]{%
\textbf{Historical note.}
This paper is based on the master's thesis of N.~Matevosyan,
written under the supervision of N.\,U.~Arakelian at Yerevan State University, Armenia, 1999.
The core results (Sections~1--4) were published as:
N.\,U.~Arakelian and N.~Matevosyan,
\textit{Three spheres theorem for harmonic functions},
Izv.\ Nats.\ Akad.\ Nauk Armenii Mat.\ \textbf{34} (1999), no.~3, 5--13;
English transl.\ in J.~Contemp.\ Math.\ Anal.\ \textbf{34} (1999), no.~3, 1--9.
The present expanded version includes additional uniqueness results
(Sections~5--6) that did not appear in the journal publication.
N.\,U.~Arakelian (1936--2023) was an Honored Scientist of the Republic of Armenia
and a member of the Institute of Mathematics of the Armenian Academy of Sciences.
First posted to arXiv in 2026.}
\renewcommand{\thefootnote}{\arabic{footnote}}

\section{Introduction}

Hadamard's well known Three-circle theorem represents itself a prototype of
the Two constants theorem. The latter found numerous applications in
different topics of Complex function theory, particularly in topics, related
with Phragmen-Lindel\"{o}f principle and uniqueness problems. It is natural
therefore, that several authors has attempted to find generalizations of
Hadamard's result in different directions. In particular, three-spheres
theorems for solutions of \textit{elliptic} equations has been received by
Landis [1] (in $L^{2}$-norm)and Agmon [2] (in $L^{\infty }$-norm); see also
the recent paper of Brummelhuis [3], containing a number of theorems of the
mentioned type. For the special case of \textit{harmonic} functions more
precise results has been received by Korevaar [4] and Korevaar and Meyers
[5]. Earlier this case has been considered by Solomentsev [6] , receiving a
three-spheres theorem of rather different form.

In paper [7] it has been developed some analogs of Hadamard's Three circle
theorem for harmonic functions of $n\geq 2$ variables and for three possibly
non-concentric \textit{balls} instead spheres (in weighted $L^{2}$-norm; see
Theorems 1-2 in [7]), with applications in problems on ``propagation of
smallness'' and uniqueness. In this paper we present a modification of that
theorem, which will let us receive some uniqueness results, which can be
found in Sections 5 and 6.

The main aim of this paper is to receive a direct analog of Hadamard's Three
circle theorem for harmonic functions (in weighted $L^{2}$-norm) in case of (%
$n-1$ -dimensional) non-concentric \textit{spheres} in $\mathbb{R}^{n}$.

The formulation and proof of the that result is presented in Section 4,
after necessary preparations in Section 2 and 3.

\section{Preliminary constructions}

For a proper subdomain $\Omega $ of $\mathbb{R}^{n}$ $(n\geq 2)$ we set $%
d(x)=dist(x,\partial \Omega )$. We associate with $\Omega $ a set $\Omega
^{\star }\subset \mathbb{R}^{n+1}$, by putting 
\begin{equation*}
\Omega ^{\star }:=\{(x,r)\in \mathbb{R}^{n+1}:x\in \Omega ,\ r\in (0,d(x))\}.
\end{equation*}

Obviously, $\Omega ^{\star }$ is an open set in $\mathbb{R}^{n+1}$ and one
can see, that it is also connected, i.e. $\Omega ^{\star }$ is a domain.
Consider for this arbitrary points $(x_{1},r_{1})$, $(x_{2},r_{2})\in \Omega
^{\star }$. We join points $x_{1},x_{2}\in \Omega $ by a Jordan arc $\gamma
\subset \Omega $. There exist a number $r>0$, with $r<\min \{r_{1},r_{2}\}$,
such that $d(x)>r$ for all $x\in \gamma $. Then $\gamma \times \{r\}\subset
\Omega ^{\star }$ and the union of connected sets $\gamma \times \{r\}$, $%
\{x_{1}\}\times \lbrack r,r_{1}]$, $\{x_{2}\}\times \lbrack r,r_{2}]$ is a
connected subset of $\Omega ^{\star }$, joining points $(x_{1},r_{1})$ and $%
(x_{2},r_{2})$.

With a function $f\in C(\Omega )$ associate a function $a\in C^{1}(\Omega
^{\star })$ by formula 
\begin{equation}
a(x,r)=a(x,r,f):=\int_{B_{x,r}}f(y)dy,\quad (x,r)\in \Omega ^{\star }, 
\tag{1}
\end{equation}
where $B_{x,r}$ is the open ball in $\mathbb{R}^{n}$ with center $x$ and
radius $r$. Further we will denote $B_{r}=B_{0,r},\ B=B^{n}:=B_{1},$ and
also $S_{x,r}:=\partial B_{x,r}$, $S_{r}:=S_{0,r}$; put especially $%
S=S^{n-1}:=S_{1}$.

For a point $x\in \mathbb{R}^{n}$ denote by $x_{k}$ its $k^{th}$ coordinate
with respect the standard basis $\{e_{k}\}_{k=1}^{n}$ of $\mathbb{R}^{n}$,
and for $1\leq k\leq n$ consider in $\mathbb{R}^{n}$ projection operators $%
p_{k}x:x\rightarrow x-x_{k}e_{k}$, so that $p_{k}x$ is independent of the $%
k^{th}$ coordinate. Let $e$ be a unit vector in $\mathbb{R}^{n}$, i.e. $e\in
S^{n-1}$, and let $\partial _{e}$ denotes the derivative in variable $x$ in
direction $e$. Put for simplicity $\partial _{k}=\partial _{e_{k}}$ - the
partial derivative with respect the $k^{th}$ coordinate variable. Writing
(1) as 
\begin{equation*}
a(x,r)=\int_{B_{p_{k}x,r}}f(y+x_{k}e_{k})dy=\int_{B_{p_{k}x,r}}\partial
_{k}\int_{x_{k}}^{x_{k}+y_{k}}f(p_{k}y+te_{k})dtdy,
\end{equation*}
we receive according to familiar divergence theorem, that

\begin{equation*}
a(x,r)=\int_{S_{p_{k}x,r}}n_{k}%
\int_{x_{k}}^{x_{k}+y_{k}}f(p_{k}y+te_{k})dtds.
\end{equation*}
Here $n_{k}$ is the $k^{th}$ coordinate of the outward unit normal $n$ and $%
s $ denotes the Lebesgue surface-area measure on $S$$_{p_{k}x,r}$ (and
further-on any sphere in $\mathbb{R}^{n}$). Using the independence of $%
S_{p_{k}x,r}$ from $k^{th}$ coordinate variable $x_{k}$, we receive

\begin{eqnarray*}
(\partial _{k}a)(x,r)
&=&\int_{S_{p_{k}x,r}}f(y+x_{k}e_{k})n_{k}ds-%
\int_{S_{p_{k}x,r}}f(p_{k}y+x_{k}e_{k})n_{k}ds= \\
&&\int_{S_{x,r}}f(y)n_{k}ds,
\end{eqnarray*}
since the second integral in middle expression equals to zero ( integrand
has the same absolute value with alternating signs at the points of $%
S_{p_{k}x,r},$ symmetric with respect the hyperline $y_{k}=0$). Thus we
arrive at the formula 
\begin{equation}
(\partial _{e}a)(x,r,f)=\int_{S_{x,r}}f(y)e\cdot nds,(x,r)\in \Omega ^{\star
}.  \tag{2}
\end{equation}
Let us denote by $\partial _{0}$ the derivative with respect the variable $%
r. $ Then formula 
\begin{equation}
(\partial _{0}a)(x,r)=\int_{S_{x,r}}fds,\quad (x,r)\in \Omega ^{\star } 
\tag{3}
\end{equation}
follows directly from the identity

\begin{equation*}
a(x,r)=\int_{0}^{r}\int_{S_{x,t}}fdsdt,\quad (x,r)\in \Omega ^{\star }.
\end{equation*}

It is immediately from (2) and (3) we can receive the inequality

\begin{equation*}
\left| (\nabla _{x}a)(x,r,f)\right| \leq (\partial _{0}a)(x,r,\left|
f\right| ),\quad (x,r)\in \Omega ^{\star },
\end{equation*}
noted in [8] for the case $f\geq 0$ (see Lemma 3 in [8]).

Suppose now $\gamma $ be a smooth Jordan arc in $\Omega $, $r\in
C^{1}(\gamma )$ and $0<r(x)<d(x)$ for $x\in \gamma .$ Subject $\gamma $ to
length parametrization on $[0,T]$, where $T$ is the length of $\gamma ,$ so
that $\mid \gamma \prime (t)\mid \equiv 1$ on $[0,T]$. Then $e=\gamma
^{\prime }(t)$ denotes the unit vector in direction of the tangent to $%
\gamma $ at a point $x=\gamma (t).$ We will also denote by $r^{\prime }$ the
tangential derivative of $r$ on $\gamma $, so that $r^{\prime }(x)=(r\bullet
\gamma )^{\prime }(t)$ for $x=\gamma (t).$

Introduce now a new function $\tilde{a}\in C^{1}(\gamma )$ , by putting 
\begin{equation}
\tilde{a}(x)=\tilde{a}(x,f):=a(x,r(x),f)=\int_{B_{x,r(x)}}f(y)dy,\quad x\in
\gamma .  \tag{4}
\end{equation}
From (2)-(4) it follows for $x\in \gamma $ 
\begin{equation}
(\partial _{e}\tilde{a})(x)=\int_{S_{x,r(x)}}f(y)(e\cdot n+r^{\prime })ds. 
\tag{5}
\end{equation}

\section{Correlated non-concentric spheres}

Specify now $\Omega =B=B^{n}$ and suppose $f\in C(\overline{B}).$ Then
formula (5) is valid also under the assumption $0<r(x)\leq d(x).$

\textsf{\textbf{\textit{Definition 1}}}. Let us say two balls $%
B_{x,r}\subset B$ and $B_{\overline{x},\overline{r}}\subset B$ are \textit{%
correlated} (over $B$), if vectors $x$ and $\overline{x}$ are codirected
(i.e. either $\left| x\right| \left| \overline{x}\right| =0$ or $x=\lambda 
\overline{x}$ with $\lambda >0$) and the correlation 
\begin{equation}
(1+\left| x\right| ^{2}-r^{2})\left| \overline{x}\right| =\left| x\right|
(1+\left| \overline{x}\right| -\overline{r}^{2})  \tag{6}
\end{equation}
be satisfied. Then corresponding spheres $S_{x,r},$ $S_{\overline{x},%
\overline{r}}$ also will be called correlated (over $S$ ). Obviously, the
condition is satisfied, if $x=\overline{x}=0,$ i.e. balls, concentric with $%
B $ are correlated.

Supposing $\left| \overline{x}\right| \leq \left| x\right| $ and taking into
account codirectness of vectors $x $ and $\overline{x}, $ one can rewrite
(6) as

\begin{equation*}
\left| x-\overline{x}\right| [1-(\left| x\right| +r)^{2}+\left| x\right|
(2r+\left| x-\overline{x}\right| )]=\left| x\right| (\overline{r}^{2}-r^{2}),
\end{equation*}
from which it follows first, that $\overline{r}\geq r.$ Further we have the
inequality

\begin{equation*}
\left| x-\overline{x}\right| (2r+\left| x-\overline{x}\right| )\leq (%
\overline{r}-r)(\overline{r}+r),
\end{equation*}
which implies $\left| x-\overline{x}\right| \leq \overline{r}-r$ (the
converse assumption will lead to contradiction!). The last inequality means
simply, that $B_{x,r}\subset B_{\overline{x},\overline{r}}(\subset B)$, so
that this inclusion stands equivalent the condition $\left| \overline{x}%
\right| \leq \left| x\right| .$

Note, that the equality $\left| x-\overline{x}\right| =\left| x\right|
-\left| \overline{x}\right| =\overline{r}-r $ holds only if $r=1-\left|
x\right| $ and then $\overline{r}=1-\left| \overline{x}\right| , $ i. e. if
the balls $B_{x,r} $ and $B_{\overline{x},\overline{r}} $ will touch $B $ in
some point. If exclude this case, then the strong inequality $\left| 
\overline{x}\right| <\left| x\right| $ implies $\overline{B}_{x,r}\subset B_{%
\overline{x},\overline{r}} $ and the additional assumption $\overline{x}\neq
0 $ implies $\overline{B}_{\overline{x},\overline{r}}\subset B. $

We will consider a three-spheres theorem for harmonic functions for the case
of non-concentric, non-touching, but correlated spheres. For this it will be
convenient for us to fix the ball $B_{x,r}$ with $x\neq 0$ and $0<r<1-\left|
x\right| ,$ guaranteeing so the balls $B_{x,r}$ and $B$ be non-concentric
and non-touching. To consider the family of balls, correlated with $B_{x,r}$
and containing the latter, we subject the interval $[0,x]\subset \mathbb{R}%
^{n}$ to length parametrization, putting $\overline{x}=$$x_{t}=\gamma (t)=te$
for $0\leq t\leq \left| x\right| ,$ and $e=\left| x\right| ^{-1}x$ . The
radius $r_{t}=\overline{r}$ we can define from (6), putting there $\left| 
\overline{x}\right| =t$ and considering it as an equation. Thus the balls $%
B_{x_{t},r_{t}}$ , $t\in \lbrack 0,\left| x\right| ],$ generated by $B_{x,r}$
via correlation, are well defined, coincided with $B_{x,r}$ for $t=\left|
x\right| $ and with $B$ for $t=0.$

Consider now a sphere $S_{a,\rho }$ and the \textit{inversion} $\varphi $\
with respect to $S_{a,\rho }$: 
\begin{equation}
\varphi (y)=a+\frac{\rho ^{2}}{\left| y-a\right| ^{2}}(y-a),  \tag{7}
\end{equation}
realizing one-to-one and conform transformation of $\mathbb{R}^{n}\backslash
\{a\}$ onto $\mathbb{R}^{n}\backslash \{a\}$ and preserving the family of
spheres in $\mathbb{R}^{n}\backslash \{a\}$ (see [10], p.60). If we choose $%
a=\left| a\right| e$ with $\left| a\right| >1$ and put $\rho ^{2}=\left|
a\right| ^{2}-1,$then $S_{a,\rho }$will be orthogonal to $S$, so that $%
\varphi (B)=B$ with $\varphi (0)=0$. It is important for our further
consideration, that $\left| a\right| $\ may be chosen in such a way, that $%
\varphi (B_{x_{t},r_{t}})=B_{r_{t}^{\ast }}$ for all $t\in \lbrack 0,\left|
x\right| ]$, i.e. the images of balls $B_{x_{t},r_{t}}$\ be \textit{%
concentric,} with centers at origin. Actually, one can find $\left| a\right|
>1$ from the quadratic equation (see equivalent formula (21) in [7])

\begin{equation}
\frac{\left| a\right| ^{2}+1}{\left| a\right| }=\frac{1+\left| x\right|
^{2}-r^{2}}{\left| x\right| }.  \tag{8}
\end{equation}

From correlation (6) (with $\left| \overline{x}\right| =t$ ) one can write $%
r_{t}$ in terms of $a$, using (8): 
\begin{equation}
r_{t}^{2}=(1-\left| a\right| ^{-1}t)(1-\left| a\right| t).  \tag{9}
\end{equation}

Also (see formula (18) in [7]) the image radius $r_{t}^{\star }$ may be
defined by formula 
\begin{equation}
r_{t}^{\star }=r_{t}\frac{\left| a\right| }{\left| a\right| -t}\equiv \frac{%
1-\left| a\right| t}{r_{t}},  \tag{10}
\end{equation}

which with (9) gives us 
\begin{equation}
(r_{t}^{\star })^{2}=\left| a\right| \frac{1-\left| a\right| t}{\left|
a\right| -t}.  \tag{11}
\end{equation}

Differentiating (9) and (11) with respect the variable $t$ and taking (10)
into account one more time , we receive 
\begin{equation}
2r_{t}r_{t}^{\prime }=2t-\frac{\left| a\right| ^{2}+1}{\left| a\right| }%
,\quad 2r_{t}(r_{t}^{\star })^{\prime }\left| a\right| =-\rho ^{2}\left( 
\frac{r_{t}^{\star }}{r_{t}}\right) .  \tag{12}
\end{equation}

We wish to apply formula (5) to our specified case, replacing previously
there $x $ by $\overline{x} $, choosing $\gamma (t)=te $ for $t\in [0,\left|
x\right| ] $ with $e=x/\left| x\right| $ (so that $\gamma ^{\prime
}(t)\equiv e $) and $r(\overline{x})=r_{t} $ for $\overline{x}=x_{t}=\gamma
(t). $ If now $y\in S_{_{x_{t},r_{t}}} $ , then $n=n_{y}=r_{t}^{-1}(y-te), $
so that from the first equality in (12) and the equality $a=\left| a\right|
e $ it follows

\begin{equation*}
e\cdot n+r_{t}^{\prime }=\frac{e\cdot y-t+r_{t}r_{t}^{\prime }}{r_{t}}=-%
\frac{\left| y-a\right| ^{2}+1-\left| y\right| ^{2}}{2\left| a\right| r_{t}}.
\end{equation*}

Substituting this into (5) and noting, that $\partial _{e}=d/dt,$ we arrive
at to formula 
\begin{equation}
\frac{d}{dt}\int_{B_{x_{t},r_{t}}}f(y)dy=-\frac{1}{2\left| a\right| r_{t}}%
\int_{S_{x_{t},r_{t}}}f(y)[\left| y-a\right| ^{2}+1-\left| y\right| ^{2}]ds,
\tag{13}
\end{equation}
valid for any function $f\in C(\overline{B})$ , if $t\in \lbrack 0,\left|
x\right| ]$.

Further we will require some estimations related to $r_{t}^{\star }. $ Put

\begin{equation*}
(r_{t}^{\star })^{2}=1-\tau ,\quad \tau =\frac{\rho ^{2}t}{\left| a\right| -t%
}
\end{equation*}
as seen from (11). It follows from (10) $\left| a\right| t<1$ for $t\in
\lbrack 0,\left| x\right| ],$ which implies $0<\tau <1.$ Applying now the
inequality $\log (1-\tau )^{-1}>\tau ,$ we receive 
\begin{equation}
\log (r_{t}^{\star })>\frac{\tau }{2}>t\frac{\left| a\right| -1}{\left|
a\right| -t}.  \tag{14}
\end{equation}

Noting, that by (8), $\left| a\right| -1>1-\left| x\right| -r$ (add $-2$ to
both sides of (8) ) and $\left| a\right| -t\leq \left| a\right| <\left|
x\right| ^{-1}\leq t^{-1}$ , we receive from (14), that 
\begin{equation}
\log (r_{t}^{\star })^{-1}>t^{2}(1-\left| x\right| -r).  \tag{15}
\end{equation}

In converse direction, it follows from (10), that $(r_{t}^{\star
})^{-1}<r_{t}^{-1}(1-\left| x\right| t).$ Taking here $t=\left| x\right| $
and using (15), we finally receive

\begin{equation}
\alpha _{t}:=\frac{\log r_{t}^{\star }}{\log r_{\left| x\right| }^{\star }}%
>\omega _{t}:=\frac{t^{2}(1-\left| x\right| -r)}{\log \frac{1-\left|
x\right| ^{2}}{r}}.  \tag{16}
\end{equation}

\section{Log-Convexity}

We set for a function $f\in C(B_{x,\tau })$ and $0<r<\tau :$ 
\begin{equation}
L_{p}(x,r,f)=\left( \int_{S_{x,r}}\left| f\right| ^{p}ds\right) ^{1/p},\quad
A_{p}(x,r,f)=\left( \int_{B_{x,r}}\left| f(y)\right| ^{p}dy\right) ^{1/p}, 
\tag{17}
\end{equation}
if $0<p<\infty $ and

\begin{equation*}
L_{\infty }(x,r,f)=\sup_{S_{x,r}}|f|,\quad A_{\infty
}(x,r,f)=\sup_{B_{x,r}}\left| f\right|
\end{equation*}
for $p=\infty .$ We will omit $x$ in this notations, if $x=0.$

\textbf{Definition2}. \textsf{A function }$L>0$\textsf{\ on }$(\tau
_{0},\tau )\subset R^{+}$\textsf{\ is said to be} \textsf{\textit{%
logarithmically}} \textsf{\textit{convex}} \textsf{(log-convex) of the
logarithm (of log) if for }$r\in \lbrack r_{1},r_{2}]\subset (\tau _{0},\tau
)$ 
\begin{equation}
L(r)\leq L^{\alpha }(r_{1})L^{1-\alpha }(r_{2}),\quad \alpha =\frac{\log
(r_{2}/r)}{\log (r_{2}/r_{1})}.  \tag{18}
\end{equation}

For $\gamma >0, $ it follows $L^{\gamma } $ is log -convex of log, provided
it is true for $L. $

\textsf{\textbf{\textit{Remark 1}}}. If log-convex function $L $ of log is
increasing, then one can replace in (18) $\alpha $ by any $\beta \in
(0,\alpha ]. $

\textsf{\textbf{\textit{Examples}}}. \textsf{\textbf{\textit{a}}}). If $f$
is a holomorphic function in $B$$^{2}\subset \mathbb{R}^{2}$, Hardy's
Convexity theorem ( which is derived from Hadamard's classical Three-circle
theorem, containing the latter for $p=\infty $), states (see [9], p.9,
Theorem 1.5) the function $L_{p}(\cdot ,f)$ for $p>0$ is log-convex of log
on $(0,1)$. The same statement is true also for the function $A_{p}(\cdot
,f) $ (see [7], Lemma 2).

\textsf{\textbf{\textit{b)}}}. The next example (important for us) related
with complex valued functions $f$, harmonic in unit ball $B^{n}\subset 
\mathbb{R}^{n},n\geq 2.$ Application of Hadamard's theorem to Parceval
identity for $f$ (see [5], Lemma 2.1) permits to state, that the
function $L_{2}(\cdot ,f)$ (and the function $A_{2}(\cdot ,f)$ as well) is
log -convex of log on $(0,1)$.

Note that in this examples the function $L $ is increasing and Remark 1 is
valid for them.

\textsf{\textbf{\textit{Remark 2}}}. In Examples \textsf{\textbf{\textit{a)}}%
} and \textsf{\textbf{\textit{b)}}} one can put in (18) $r_{2}=1$ for $%
L=L_{p}(\cdot ,f)$ and $L=L_{2}(\cdot ,f)$ respectively, assuming $f\in
H^{p}(B)$ for case \textsf{\textbf{\textit{a)}}} and $f\in h^{2}(B)$ in case 
\textsf{\textbf{\textit{b)}}}; here $H^{p}(B)$ are Hardy classes of
holomorphic functions in the unit disc $B$ (see [9]) and $h^{2}(B)$ is
the corresponding Hardy class of harmonic functions in $B=B^{n},n\geq 2$
(see [10], Ch.6). One can extend $f$ in these cases a.e. onto $S=S^{n}$
by Fatou Limit theorem and define $L_{p}(1,f)$ as in (17). Applying (18) to
corresponding $L_{p}(\cdot ,f_{\tau })$ with $0<\tau <1,$ where $f_{\tau
}(y)=f(\tau y)$, and letting $\tau \rightarrow 1,$ one can receive the
result.

The first part of the remark is valid also for the quantities $A_{p}(.,f)$
in corresponding Bergman classes $b^{p}(B)$ of holomorphic and harmonic
functions.

\section{\textbf{A THREE-SPHERES THEOREM}.}

As in Section 2, let $B_{x,r}$ be a ball in $B=B^{n},$ $n\geq 2,$
non-concentric and non-touching with $B.$ Let $a$ be a point, mentioned
there, with $\left| a\right| >1,$ codirected with $x$ and satisfying (8).
Consider a complex valued function $f,$ harmonic in a ball $B_{\tau }\subset 
\mathbb{R}^{n}$ with $1<\tau <\left| a\right| $ and denote by $f^{\star }$
the \texttt{\textit{Kelvin transform}} of $f$ with respect the inversion $%
\varphi =\varphi ^{-1},$ defined in (7): 
\begin{equation}
f^{\star }=\left( \frac{\rho }{\left| y-a\right| }\right) ^{n-2}\overline{%
f(\varphi (y))},\quad y\in \varphi (B_{\tau }).  \tag{19}
\end{equation}

Since $\varphi (B)=B,$ there exists a $\tau _{1}>1,$ such that $B_{\tau
_{1}}\subset \varphi (B_{\tau }).$ Note, that $f^{\star }$ is harmonic in $%
\varphi (B_{\tau })$ and particularly in $B_{\tau _{1}};$ if $n=2$ and $f$
is holomorphic function in $B_{\tau },$ then $f^{\star }$ is holomorphic in $%
\varphi (B_{\tau }).$

Consider now the family $\left\{ B_{x_{t},r_{t}}\right\} $ of balls for $%
t\in \lbrack 0,\left| x\right| ],$ correlated with $B_{x,r}$ (over $B$ ) and
containing $B_{x,r}.$ We apply to function $L_{2}^{2}(\cdot ,f^{\star })$
inequality (18) for arguments $r_{\left| x\right| }^{\star }\leq
r_{t}^{\star }\leq 1,$ replacing there $\alpha $ by any $\beta \in (0,\alpha
]$ (see Example 2 and Remark 1). Reminding, that by (18), 
\begin{equation}
r_{t}^{\star }=\left( r_{\left| x\right| }^{\star }\right) ^{\alpha }\leq
\left( r_{\left| x\right| }^{\star }\right) ^{\beta },  \tag{20}
\end{equation}
we receive 
\begin{equation}
L_{2}^{2}\left( r_{t}^{\star },f^{\star }\right) \leq L_{2}^{2\beta }\left(
r_{\left| x\right| }^{\star },f^{\star }\right) L_{2}^{2(1-\beta )}\left(
1,f^{\star }\right) .  \tag{21}
\end{equation}

We note now, in view of (19), that by formula for the change of variables ($%
\zeta =\varphi (y) $) in integrals with volume measures,

\begin{equation*}
I(t):=\int_{B_{r_{t}}^{\star }}\left| f^{\star }(y)\right|
^{2}dy=\int_{B_{x_{t},r_{t}}}\rho ^{4}\left| \zeta -a\right| ^{-4}\left|
f(\zeta )\right| ^{2}d\zeta ,
\end{equation*}
taking into account, that $\varphi (B_{x_{t},r_{t}})=B_{r_{t}^{\star }}$ and 
$\det \varphi ^{\prime }(y)=-\left| y-a\right| ^{-2n}.$

To express (21) directly in terms of the function $f$ and spheres $%
S_{x_{t},r_{t}},$ $S_{x,r}$ and $S,$ we first apply formula (3) and then
formula (13). Using also the second formula in (12), we receive 
\begin{equation}
L_{2}^{2}\left( r_{t}^{\star },f^{\star }\right) =\frac{1}{\left(
r_{t}^{\star }\right) ^{\prime }}I^{\prime }(t)=\rho ^{2}\frac{r_{t}}{%
r_{t}^{\star }}\int_{S_{x_{t},r_{t}}}\left| f(y)\right| ^{2}ds_{a},  \tag{22}
\end{equation}
where the measure $s_{a}$ related with $s$ by formula 
\begin{equation}
ds_{a}=\frac{\left| y-a\right| ^{2}+1-\left| y\right| ^{2}}{\left|
y-a\right| ^{4}}ds.  \tag{23}
\end{equation}

Let us substitute (23) into (20), taking into account (21), and come back
again to notations in the beginning of Section 2: $x_{t}=\overline{x},r_{t}=%
\overline{r}.$ We finally arrive at the inequality 
\begin{equation}
\overline{r}\int_{S_{\overline{x},\overline{r}}}\left| f\right|
^{2}ds_{a}\leq \left( r\int_{S_{x,r}}\left| f\right| ^{2}ds_{a}\right)
^{\beta }\left( \int_{S}\left| f\right| ^{2}ds_{a}\right) ^{1-\beta }, 
\tag{24}
\end{equation}
valid for any $\beta \in (0,\alpha ],$ where $(\overline{r})^{\star
}=(r^{\star })^{\alpha },$ or by (10), 
\begin{equation}
\frac{1-\left| a\right| \left| \overline{x}\right| }{\overline{r}}=\left( 
\frac{1-\left| a\right| \left| x\right| }{r}\right) ^{\alpha }.  \tag{25}
\end{equation}
In particular, it follows from the estimate (16), that one can put in (24)

\begin{equation*}
\beta =\omega :=\left| \overline{x}\right| ^{2}\frac{1-\left| x\right| -r}{%
\log \frac{1-\left| x\right| ^{2}}{r}}.
\end{equation*}

Note also one can easily extend the inequality (24) for functions $f $ from
the Hardy space $h^{2}(B^{n}) $ of harmonic functions in the way, mentioned
in Remark 2.

We summarize our discussion in the following theorem.

\texttt{\textbf{\textit{Theorem 1.}}} \textit{Let} $S_{x,r}$\textit{\ and }$%
S_{\overline{x},\overline{r}}$ \textit{be correlated spheres with respect to 
}$S=S^{n-1}$ $(n\geq 2)$\textit{\ with }$x\neq 0,$\textit{\ }$0<r<1-\left|
x\right| $\textit{\ and }$\left| \overline{x}\right| \leq \left| x\right| $%
\textit{$.$ Then for a function }$f\in h^{2}(B^{n})$\textit{\ inequality
(24) holds for any} $\beta \in (0,\alpha ]$\textit{\ with }$\alpha ,$ 
\textit{satisfying (25), and measure} $s_{a}$\textit{\ is defined by (23);
the point }$a=\left| a\right| \left| x\right| ^{-1}x$\textit{\ with }$\left|
a\right| >1$\textit{\ may be found from (8).}

\texttt{\textbf{\textit{Remark~3}}}. In case $n=2,$ assuming $f$ be
holomorphic and actually $f\in H^{2}(B^{2}),$ in (24) $ds_{a}$ can be
replaced $(\left| y-a\right| ^{2}+1-\left| y\right| ^{2})ds.$ This follows
from (24) by replacing $f(y)$ to $(y-a)^{2}f(y).$

\bigskip

\section{Three balls theorem (non-concentric case)}

Earlier an analog of Theorem 1 but with balls instead of spheres was proved
in [7]. That theorem is the following:

\textbf{Theorem 2}\texttt{.} \textit{Let }$u\in h^{2}(B)$\textit{, }$%
B_{x_{0},r_{0}}$\textit{\ and }$B_{\overline{x},\overline{r}}$\textit{\ be
correlated balls with respect to }$B=B^{n}$\textit{\ }$(n\geq 2)$\textit{\ .
Denote }

\begin{equation}
\delta _{0}:=\frac{\left| \overline{x}\right| ^{2}}{2(1-\left| \overline{x}%
\right| )}\frac{1-\left| x_{0}\right| -r_{0}}{\log \frac{1-\left|
x_{0}\right| ^{2}}{r_{0}/2}}.  \tag{26}
\end{equation}
\textit{\ }

\textit{Then for every }$\delta \in (0,\delta _{0}]$\textit{\ the following
holds }

\begin{equation}
\int_{B_{\overline{x},\overline{r}}}\left| u\right| ^{2}d\mu _{a}\leq \left(
\int_{B_{x_{0},r_{0}}}\left| u\right| ^{2}d\mu _{a}\right) ^{\delta }\left(
\int_{B}\left| u\right| ^{2}d\mu _{a}\right) ^{1-\delta }\mathit{,}  \tag{27}
\end{equation}

\textit{where }

\begin{equation}
d\mu _{a}=\left| x-a\right| ^{-4}dx\mathit{.}  \tag{28}
\end{equation}

\bigskip

For some applications it is useful to have (27) in terms of $dx$ rather than 
$d\mu _{a}$. That is done in the following theorem by using the technique of
Imbedding the space $\mathbb{R}^{n}$ into $\mathbb{R}^{n+5},$ used in proof
of Theorem 6 in [7].

\textbf{Statement 1}\texttt{.} \textit{Let }$u\in h^{2}(B)$\textit{, }$%
\left| x_{0}\right| \geq 1/2,$ $B_{x_{0},r_{0}}$\textit{\ and }$B_{\overline{%
x},\overline{r}}$\textit{\ be correlated balls with respect to }$B=B^{n}$%
\textit{\ }$(n\geq 2),$ \textit{and }$\delta _{0}$\textit{\ is as in (26).
Then for every }$\lambda \in (0,1)$ \textit{and }$\delta \in (0,\delta
_{0}], $\textit{\ the following holds}

\begin{equation}
\int_{B_{\overline{x},\lambda \overline{r}}}\left| u\right| ^{2}dx\leq \frac{%
405}{\overline{r}(1-\lambda ^{2})^{5/2}}\left( \int_{B_{x_{0},r_{0}}}\left|
u\right| ^{2}dx\right) ^{\delta }\left( \int_{B}\left| u\right|
^{2}dx\right) ^{1-\delta }\mathit{.}  \tag{29}
\end{equation}

\bigskip

For the proof let us imbed the space $\mathbb{R}^{n}$ into $\mathbb{R}^{n+5}$
by setting the points of $\mathbb{R}^{n+5}$ by $(x,y),$ where $x\in \mathbb{R%
}^{n},$ $y\in \mathbb{R}^{5}.$ Let us continue $u\in h^{2}(B^{n})$ to the
function $\tilde{u}\in h^{2}(B^{n+5})$ in the following way

\begin{equation*}
\tilde{u}(x,y)=u(x),\quad x\in B^{n},\quad (x,y)\in B^{n+5}.
\end{equation*}

Let $\nu _{n}$ denote the $n$-dimensional Lebesgue measure, 
\begin{equation*}
\nu _{n}(r)=\nu (B_{x,r})=\nu _{n}r^{n}.
\end{equation*}
For the function $\tilde{u}$ we have

\begin{eqnarray*}
\int_{B^{n+5}}\left| \tilde{u}\right| ^{2}dxdy &=&\int_{B^{n+5}}\left|
u\right| ^{2}dxdy=\int_{B^{n}}\left| u\right| ^{2}\nu _{5}(B_{\sqrt{1-\left|
x\right| ^{2}}}^{5})dx\leq \\
&\leq &\nu _{5}(B^{5})\int_{B^{n}}\left| u\right| ^{2}dx<\infty .
\end{eqnarray*}

We are imbedding the ball $B_{x,r}^{n}$ into the ball $%
B_{(x_{0},0_{5}),r_{0}}^{n+5},$ the ball $B_{\bar{x},\bar{r}}^{n}$ into the
ball $B_{(\bar{x},0_{5}),\bar{r}}^{n+5}$ and the ball $B^{n}$ into the ball $%
B^{n+5},$ where $0_{5}=(0,0,0,0,0).$ We get 
\begin{equation*}
B_{(x_{0},0_{5}),r_{0}}^{n+5}\subset B_{(\bar{x},0_{5}),\bar{r}%
}^{n+5}\subset B^{n+5}.
\end{equation*}
Instead of $a$ in $\mathbb{R}^{n+5}$ we will consider $(a,0_{5}).$ For the
points $x\in \mathbb{R}^{n},$ $y\in \mathbb{R}^{5}$ we get

\begin{equation*}
d\mu _{(a,0_{5})}=\left| (x,y)-(a,0_{5})\right| ^{-4}dxdy=\left| \left|
x-a\right| ^{2}+\left| y\right| ^{2}\right| ^{-2}dxdy\mathit{.}
\end{equation*}

Observe that for every $b\in \mathbb{R}^{n},$ $l\in (0,1]$\ and $g(x,y)\in
C(B_{(b,0_{5}),l}^{n+5})$\ it is true, that

\begin{equation}
\int_{B_{(b,0_{5}),l}^{n+5}}g(x,y)dxdy=\int_{B_{b,l}^{n}}\left( \int_{B_{%
\sqrt{l-\left| x-b\right| ^{2}}}^{5}}g(x,y)dy\right) dx=  \tag{30}
\end{equation}
\begin{equation*}
=\int_{B_{b,l}^{n}}\left( \int_{0}^{\sqrt{l-\left| x-b\right| ^{2}}%
}t^{4}\left( \int_{S}g(x,t\theta )d\theta \right) dt\right) dx
\end{equation*}
\ \ 

\bigskip The following estimates will be useful in the future:

\begin{equation}
\left( \left| x-a\right| ^{2}+\left| y\right| ^{2}\right) ^{-2}\leq \left|
y\right| ^{-4}  \tag{31}
\end{equation}

and if $y\in B_{\sqrt{\bar{r}^{2}-\left| x-\bar{x}\right| ^{2}}}^{5},$ then

\begin{equation}
\left( \left| x-a\right| ^{2}+\bar{r}^{2}-\left| x-\bar{x}\right|
^{2}\right) ^{-2}\leq \left( \left| x-a\right| ^{2}+\left| y\right|
^{2}\right) ^{-2}.  \tag{32}
\end{equation}

Applying the first equation of (30) for the left hand side and the whole
(30) for the right hand side of (27), we will receive

\begin{gather}
\frac{1}{5}\int_{B_{\overline{x},\overline{r}}^{n}}\left| \tilde{u}%
(x,0_{5})\right| ^{2}\frac{\left( \bar{r}^{2}-\left| x-\bar{x}\right|
^{2}\right) ^{5/2}}{\left( \left| x-a\right| ^{2}+\bar{r}^{2}-\left| x-\bar{x%
}\right| ^{2}\right) ^{2}}dx\leq  \tag{33} \\
\leq \left( \int_{B_{x_{0},r_{0}}^{n}}\left| \tilde{u}(x,0_{5})\right|
^{2}\left( r_{0}^{2}-\left| x-x_{0}\right| ^{2}\right) ^{1/2}dx\right)
^{\delta }  \notag \\
\cdot \left( \int_{B}\left| \tilde{u}(x,0_{5})\right| ^{2}\left( 1-\left|
x\right| ^{2}\right) ^{1/2}dx\right) ^{1-\delta }.
\end{gather}

To estimate the right part of (33), let us note, that

\begin{equation*}
\left( r_{0}^{2}-\left| x-x_{0}\right| ^{2}\right) ^{1/2}\leq 1\qquad \left(
1-\left| x\right| ^{2}\right) ^{1/2}\leq 1.
\end{equation*}

to estimate the left part of (33), let us use the following estimate, which
has been noted in Theorem 5 of [7]:

\begin{equation*}
\left| (x,y)-(a,0_{5})\right| ^{2}=\left| x-a\right| ^{2}+\left| y\right|
^{2}\leq 3^{2}\quad \text{if }\left| x_{0}\right| \geq 1/2.
\end{equation*}

From here we get

\begin{equation*}
\left( \left| x-a\right| ^{2}+\bar{r}^{2}-\left| x-\bar{x}\right|
^{2}\right) ^{-2}\geq 81^{-1}.
\end{equation*}

Now let us take $\lambda \in (0,1)$ and consider $B_{\overline{x},\lambda 
\overline{r}}\subset B_{\overline{x},\overline{r}}$ instead of $B_{\overline{%
x},\overline{r}}$ in the left hand side of (33). For $x\in B_{\overline{x}%
,\lambda \overline{r}}$, it is true, that $\left| x-\bar{x}\right| <\lambda 
\bar{r}$. From all these estimates we receive (29). The proof is complete.

Simple calculations can give us the analog of Theorem 3,for the case, when
instead of $B=B^{n}$ we consider a ball $B_{R}=B_{R}^{n}$ with the centre at
the origin but arbitrary radius $R$. We will have the following :

\textbf{Statement 2}\texttt{.} \textit{Let }$u\in h^{2}(B_{R})$\textit{, }$%
\left| x_{0}\right| \geq 1/2,$\textit{\ }$B_{x_{0},r_{0}}\subset B_{R}.$ 
\textit{Let }$\bar{x}\in \lbrack 0,x_{0}],$ $\bar{r}$\textit{\ be root of
the following equation }

\begin{equation}
(R^{2}+\left| x_{0}\right| ^{2}-r_{0}^{2})\left| \overline{x}\right| =\left|
x_{0}\right| (R^{2}+\left| \overline{x}\right| -\overline{r}^{2})  \tag{34}
\end{equation}

\textit{\ and }$\delta _{0}$\textit{\ is defined by }

\begin{equation}
\delta _{0}:=\frac{\left| \overline{x}\right| ^{2}}{2R^{2}(R-\left| 
\overline{x}\right| )}\frac{R-\left| x_{0}\right| -r_{0}}{\log \frac{%
R-\left| x_{0}\right| ^{2}}{r_{0}/2}}.  \tag{35}
\end{equation}

\textit{Then for every }$\lambda \in (0,1)$ \textit{and }$\delta \in
(0,\delta _{0}],$\textit{\ the following holds}

\begin{equation}
\int_{B_{\overline{x},\lambda \overline{r}}}\left| u\right| ^{2}dx\leq \frac{%
405}{(1-\lambda ^{2})^{5/2}}\left( \frac{R}{\overline{r}}\right) ^{5}\left(
\int_{B_{x_{0},r_{0}}}\left| u\right| ^{2}dx\right) ^{\delta }\left(
\int_{B_{R}}\left| u\right| ^{2}dx\right) ^{1-\delta }\mathit{.}  \tag{36}
\end{equation}

\bigskip

Let us make some notations:

$
\begin{array}{ll}
A_{p}(x,r,f):=\left( \frac{1}{\nu (B_{x,r})}\int_{B_{x,r}}\left| f\right|
^{p}d\nu \right) ^{1/p} & 0<p<\infty , \\ 
A_{p}(r,f):=A_{p}(0,r,f) & 0<p<\infty .
\end{array}
$

\texttt{\textbf{\textit{Remark~3}}}. Rewriting (36) in terms of $%
A_{2}(x,r,u) $ we receive the following inequation:

\begin{equation}
A_{2}(\bar{x},\lambda \bar{r},u)\leq \frac{\sqrt{405}}{(1-\lambda ^{2})^{5/4}%
}\left( \frac{R}{\overline{r}}\right) ^{\frac{n+5}{2}}A_{2}^{\delta
}(x_{0},r_{0},u)A_{2}^{1-\delta }(R,u)  \tag{37}
\end{equation}

\bigskip

Now we can receive some uniqueness properties for harmonic in $\mathbb{R}^{n}
$ functions.

\bigskip

\textbf{Theorem 3.}\textit{\ Let }$u\in h(\mathbb{R}^{n})$\textit{\ }$(n\geq
2),$\textit{\ }$A_{2}(r,u)\leq \varphi (r),$\textit{\ where }$\varphi $%
\textit{\ is a monotonic increasing function in }$[0,\infty ).$\textit{\
Assume that }$A_{2}(x_{m},r_{m},u)\leq \varepsilon _{m}$\textit{, when }$%
m\in N$\textit{\ and }$0<2r_{m}\leq \left| x_{m}\right| .$\textit{\ Then if} 
\begin{equation}
\lim_{m\rightarrow \infty }\left[ \frac{\rho _{m}}{100}\log \varepsilon
_{m}+\varphi (4\left| x_{m}\right| )\right] =-\infty \mathit{,\quad }\text{%
\textit{where }}\rho _{m}=\frac{1}{\log \frac{2\left| x_{m}\right| }{r_{m}}},
\tag{38}
\end{equation}

\textit{then }$u\equiv 0$\textit{\ in }$R$\textit{\ }$^{n}.$

\bigskip

The proof follows using contradictory argument. So suppose the conclusion in
theorem fails. Take $\forall r\in \lbrack 0,+\infty )$ and $\forall \eta \in
(0,+\infty ).$ Using compactness argument, we can receive, that there exists
a number $a_{r,\eta }>0$ such that

\begin{equation}
a_{r,\eta }<A_{2}(x,\eta ,u),\qquad \forall x\in \bar{B}_{r}.  \tag{39}
\end{equation}

There can be two possible cases.

\textbf{Case 1)} Suppose the sequence $\{x_{m}\}$ is bounded. Denote

\begin{equation*}
r=\overline{\lim_{m\rightarrow \infty }}\left| x_{m}\right| <+\infty .
\end{equation*}

There exists a $\eta >0$ number such that $r_{m}<\eta $ and $\eta >\frac{2r}{%
e}\geq \frac{2\left| x_{m}\right| }{e}$, when $m\in \mathbb{N}$. Consider
the following concentric balls

\begin{equation*}
B_{x_{m},r_{m}}\subset B_{x_{m},\eta }\subset B_{x_{m},e\eta }.
\end{equation*}

Since the function $A_{2}(r)=A_{2}(x,r,u)$ is log-convex of log (see [7]),
the following equation holds:

\begin{equation*}
A_{2}(x_{m},\eta ,u)\leq A_{2}^{\alpha _{m}}(x_{m},r_{m},u)A_{2}^{1-\alpha
_{m}}(x_{m},e\eta ,u),
\end{equation*}

where 
\begin{equation*}
\alpha _{m}=\frac{\log \frac{e\eta }{\eta }}{\log \frac{e\eta }{r_{m}}}=%
\frac{1}{\log \frac{e\eta }{r_{m}}}.
\end{equation*}

According to (39) there exists $a_{r,\eta }>0$ such that $A_{2}(x_{m},\eta
,u)>a_{r,\eta }$ and $b_{\eta }\geq \max \{1,\varphi (r+e\eta )\}$
satisfying $A_{2}(x_{m},\eta ,u)\leq b_{\eta }.$ From that inequation can be
received the following

\begin{equation*}
\log \frac{a_{r,\eta }}{b_{\eta }}\leq \frac{\log \varepsilon _{m}}{\log 
\frac{2\left| x_{m}\right| }{r_{m}}},
\end{equation*}

inequation, the left part of which do not depend on $m$, and the right part
is tending to $-\infty $ by (38). The received inequation is proving the
case 1).

\bigskip

\textbf{\ Case 2).} Suppose the sequence $\{x_{m}\}$ is not bounded, i.e.
there exists a subsequence $\{x_{m_{k}}\}\subset \{x_{m}\}$ such that $%
x_{m_{k}}\longrightarrow \infty $, when $k\longrightarrow \infty $. For
convenience, suppose $x_{m}\longrightarrow \infty $. Take $R_{m}=2\left|
x_{m}\right| $, $\bar{x}_{m}=3^{-1}x_{m}$ and let $\bar{r}_{m}$ to be a root
of the equation (34) taking $R=R_{m},$ $r_{0}=r_{m},$ $x_{0}=x_{m}$ and $%
\bar{x}=\bar{x}_{m}$. According to [7] we will have

\begin{gather*}
B_{x_{m},r_{m}}\subset B_{\bar{x}_{m},\bar{r}_{m}}\subset B_{R_{m}} \\
\bar{r}_{m}>\left| x-x_{m}\right| =\frac{2}{3}\left| x_{m}\right| .
\end{gather*}

Take $\delta _{m}$\ from formula (35) (taking $R=R_{m},$ $r_{0}=r_{m},$ $%
x_{0}=x_{m}$\ and $\bar{x}=\bar{x}_{m}$) we will receive 
\begin{equation*}
\delta _{m}\geq \frac{\rho _{m}}{100}.
\end{equation*}

Applying (37), for every $\lambda \in (0,1)$ we will have the following
inequality:

\begin{equation}
A_{2}(\bar{x}_{m},\lambda \bar{r}_{m},u)\leq \frac{\sqrt{405}}{(1-\lambda
^{2})^{5/4}}\left( \frac{R_{m}}{\overline{r}_{m}}\right) ^{\frac{\rho _{m}}{%
100}}A_{2}^{\delta }(x_{m},r_{m},u)A_{2}^{1-\frac{\rho _{m}}{100}}(R_{m},u).
\tag{40}
\end{equation}

Without loosing generality, we can assume, that $\varphi (4\left|
x_{m}\right| )>1,$ for every $m\in \mathbb{N}$. Take $\lambda =\frac{3}{5}.$%
\ Taking into account, that $\bar{r}_{m}>\frac{2}{3}\left| x_{m}\right| ,$\
from (40)\ we get the following inequality

\begin{equation}
A_{2}(\bar{x}_{m},\frac{3}{5}\bar{r}_{m},u)\leq 6^{\frac{n+9}{2}}\varepsilon
_{m}^{\frac{\rho _{m}}{100}}\varphi (4\left| x_{m}\right| ).  \tag{41}
\end{equation}

There exists $m_{0}\in \mathbb{N}$ such that if $m\geq m_{0},$ then $\left|
x_{m}\right| >15.$\ Observe, that $B=B_{0,1}\subset B_{\bar{x}_{m},\frac{3}{5%
}\bar{r}_{m},}$\ for $m\geq m_{0}.$\ From here, one can easily see, that

\begin{equation*}
A_{2}(1,u)\leq A_{2}(\bar{x}_{m},\frac{3}{5}\bar{r}_{m},u),\qquad \text{when 
}m\geq m_{0}.
\end{equation*}

From (41) and the last inequality one can obtain that

\begin{equation*}
6^{-\frac{n+9}{2}}A_{2}(1,u)<\varepsilon _{m}^{\frac{\rho _{m}}{100}}\varphi
(4\left| x_{m}\right| ).
\end{equation*}

Denote $a_{n}:=6^{-\frac{n+9}{2}}A_{2}(1,u).$\ We will have

\begin{equation*}
\log a_{n}<\frac{\rho _{m}}{100}\log \varepsilon _{m}+\log \varphi (4\left|
x_{m}\right| )
\end{equation*}

inequality. Since $u$\ is not identically $0$\ and $a_{n}>0$\ depends only
on $n,$ according to (38) the right side of the above inequation tends to $%
-\infty ,$ when $m\rightarrow \infty .$ That is a contradiction. The theorem
is proved. \ 

\

\end{document}